\documentclass[10pt]{amsart}
\usepackage[T1]{fontenc}
\usepackage{amsmath,amssymb,amsthm,titlesec,fancyhdr,enumerate,stackengine}
\usepackage{mathrsfs}
\usepackage{scalerel}
\usepackage{tikz}
\usepackage[margin=1in]{geometry}
\usepackage{pgfplots}
\pgfplotsset{compat=newest}
\usetikzlibrary{decorations.markings}
\usetikzlibrary{chains,shapes.multipart}
\usetikzlibrary{shapes,calc, arrows,}
\usetikzlibrary{automata,positioning,fit}
\usepackage{graphicx}	
\usepackage{caption}
\usepackage{subcaption}
\usepackage{float} 
\usepackage{listings}
\usepackage{hyperref}
\usepackage{newtxtext}
\usepackage{newtxmath}

\makeatletter
\patchcmd{\maketitle}{\@fnsymbol}{\@alph}{}{}  
\makeatother

\title{\scshape Rings of Typed Ordered Fuzzy Numbers}
	\author{M.~Kukla} 
	\thanks{\href{mailto:mkukla1@umd.edu}{mkukla1@umd.edu}}
    \author{ R.~Traylor}
    \thanks{\href{mailto:rachel.traylor@marquette.edu}
    {rachel.traylor@marquette.edu}}

	\date{\today}
 
\makeatletter
\let\runauthor\@author 
\let\runtitle\@title
\makeatother

\usepackage[backend=biber,sortcites,style=numeric,giveninits=true]{biblatex}
\DeclareNameAlias{default}{family-given}

\titleformat{\section}
  {\centering\scshape\Large}
  {\S\thesection.}
  {6pt}
  {}

\titleformat{\subsection}[runin]
  {\scshape}
  {\thesubsection.}
  {1em}
  {\addperiod}

\newtheoremstyle{solution}
{}
{}
{\color{red}}
{}
{\bfseries}
{.}
{.5em}
{}

\theoremstyle{definition}
\newtheorem{example}{Example}[section]
\newtheorem{defin}{Definition}[section]

\theoremstyle{theorem}
\newtheorem{theorem}{Theorem}[section]

\theoremstyle{solution}

\newcommand{\inbox}[1]{\mathrel{\fbox{#1}}}
\newlength\squareheight
 \newcommand\squareslash{\tikz{\draw (0,0) rectangle (\squareheight,\squareheight);\draw(0,0) -- (\squareheight,\squareheight)}}

\renewcommand{\kappa}{\varkappa}
\newcommand{\addperiod}[1]{#1.}

\begin{document}

\maketitle
\begin{abstract}
We construct rings of typed ordered fuzzy numbers whose component functions are of a common form.  As this ring also contains improper fuzzy numbers (OFNs whose membership "functions" are actually just relations), we develop a set of operations to convert an improper fuzzy number to a proper one of the same type.
\end{abstract}

\section{Introduction}
Fuzzy numbers have garnered a wide use, particularly in inference and control, but have also found statistical applications~\cite{MB2013}. A fuzzy number is a fuzzy subset $A$ in which the membership function $\mu:\mathbb{R} \to [0,1]$ satisfies the following properties~\cite{Bede2013, KaufmanGupta1985}:
\begin{enumerate}[(i)]
\item there exists $x_{0} \in \mathbb{R}$ such that $\mu(x_{0}) = 1$
\item $\mu(tx + (1-t)y) \geq \min\{\mu(x), \mu(y)\}$ for all $t \in [0,1]$, $x,y \in \mathbb{R}$
\item $\mu$ is upper semicontinuous
\item $cl(\{x \in \mathbb{R} : \mu(x) > 0\})$ is compact, where $cl(\cdot)$ denotes the closure of a set
\end{enumerate}

Dubois and Prade~\cite{DP1978} created a more restrictive class of $L-R$ fuzzy numbers, wherein the fuzzy number's membership function is created from two "spread" functions, $L,R:[0,1] \to [0,1]$ in the following way:

Let $L,R:[0,1] \to [0,1]$ be two continuous, increasing functions such that $L(0) = 0 = R(0)$, $L(1) = 1 = R(1)$. Then let $a_{0}^{-} \leq a_{1}^{-} \leq a_{1}^{+} \leq a_{0}^{+}$ all be real numbers. The fuzzy set $A$ with membership function $u : \mathbb{R} \to [0,1]$ is an \textit{L-R fuzzy number} if
\[u(x) = \begin{cases}0, & x < a_{0}^{-} \\
                  L\left(\frac{x-a_{0}^{-}}{a_{1}^{-} - a_{0}^{-}}\right), & a_{0}^{-} \leq a < a_{1}^{-} \\
                1, & a_{1}^{-} \leq x < a_{1}^{+} \\
             R\left(\frac{a_{0}^{+}-x}{a_{0}^{+} - a_{1}^{+}}\right), & a_{1}^{+} \leq x < a_{0}^{+} \\
         0, & a_{0}^{+} \leq x
\end{cases}\]

The parts of an $L-R$ fuzzy number are as follows: 

\begin{itemize}
\item the core is given by $[a_{1}^{-},a_{1}^{+}]$ 
\item the support $\text{supp}(A)$ is given by $[a_{0}^{-}, a_{0}^{+}]$  
\item  the \textit{left spread}  is $\underline{a} = a_{1}^{-} - a_{0}^{-}$, and 
\item the \textit{right spread} is $\bar{a} = a_{0}^{+}- a_{1}^{+}$ . 
\end{itemize}

For example, if $L$, $R$ are linear, then the membership function of an $L-R$ fuzzy number yields a trapezoid, and hence we call these \textit{trapezoidal fuzzy numbers}. Specifically, the membership function of a trapezoidal fuzzy number is given by 

$$u(x) = \begin{cases}0, & x < a_{0}^{-} \\
                  \frac{x-a_{0}^{-}}{a_{1}^{-} - a_{0}^{-}}, & a_{0}^{-} \leq a < a_{1}^{-} \\
                1, & a_{1}^{-} \leq x < a_{1}^{+} \\
             \frac{a_{0}^{+}-x}{a_{0}^{+} - a_{1}^{+}}, & a_{1}^{+} \leq x < a_{0}^{+} \\
         0, & a_{0}^{+} \leq x
\end{cases}$$
\subsection{Positive and Negative Fuzzy Numbers}
\label{subsec: pos neg fn}

We call a fuzzy number $A$ \textbf{positive} if all $x \in \text{supp}(A), x > 0$, and \textbf{negative} if all $x \in \text{supp}(A) < 0$. If $0 \in \text{supp}(A)$, $A$ is neither positive nor negative. 

\subsection{Arithmetic on Classical Fuzzy Numbers}
\label{subsec: arithmetic classical fn}
If fuzzy numbers are to be used in any sort of arithmetic, analytic, or data processing sense, it is necessary to define arithmetic operations on them. Arithmetic operations have been defined for $L-R$ fuzzy numbers by operating on their \textit{level sets}. The level set of an $L-R$ fuzzy number is given by

\[A_{\alpha} = [a_{0}^{-} + L^{-1}(\alpha)\underline{a}, a_{0}^{+}-R^{-1}(\alpha)], \alpha \in [0,1]\]

The sum of two fuzzy numbers $A$, $B$ is the fuzzy number $C$ whose
$\alpha-$level sets are given by 
\[C_{\alpha} = (A+B)_{\alpha} = [a_{\alpha}^{-} + b_{\alpha}^{-}, a_{\alpha}^{+} + b_{\alpha}^{+}].\] Equivalently, we can define addition on the membership functions. Let $\mu_{A}, \mu_{B}, \mu_{C}$ be the membership functions of $A,B$, and $C=A+B$. Then 
\[\mu_{C}(A)(z) = \bigvee_{z=x+y}(\mu_{A}(x) \wedge \mu_{B}(y)) \qquad \text{for all } x,y,z \in \mathbb{R}\]
where $\vee$ denotes maximum, and $\wedge$ denotes minimum~\cite{KaufmanGupta1985}.

Unfortunately, one of the biggest drawbacks to the $L-R$ fuzzy number theory and the so-defined addition operation is a lack of a unique additive inverse, which makes defining the operation of subtraction difficult~\cite{Bede2013}. Under the above criterion, $L-R$ fuzzy numbers with the addition operation as defined only produce a commutative monoidal structure, which isn't strong enough for many applications. 

Multiplication of two $L-R$ fuzzy numbers can be defined according to their membership functions or their level sets and is based on Zadeh's extension principle~\cite{Zadeh1965}. The level sets of $C = AB$ are given by 
\[C_{\alpha} = [\inf\{xy : x \in A_{\alpha}, y \in B_{\alpha}\}, \sup\{xy : x \in A_{\alpha}, y \in B_{\alpha}\}]\]

The classical number $1$ is the multiplicative identity, but unfortunately we again only have a commutative monoidal structure on the set of $L-R$ fuzzy numbers under this multiplication operation. We lack the ability to create a multiplicative inverse. 

Other issues with $L-R$ fuzzy numbers include uncontrollable results in repeatedly applied operations, which are caused by the need of intermediate approximations~\cite{Wagenknecht2001}.

\subsection{Kosinski's Ordered Fuzzy Numbers}
\label{subsec: kosinski ofn}
Kosinski~\cite{Kosinski2003, Kosinski2003b, Kosinski2006}, in an attempt to circumvent the issues described in~\ref{subsec: arithmetic classical fn}, created the notion of an \textit{ordered fuzzy number}, which resembles an $L-R$ fuzzy number, but adds an orientation. The main motivation for adding this orientation is to resolve the issue of a lack of unique additive inverse. With orientation, a true additive inverse $-A$ could be defined for a fuzzy number $A$ and addition operation $+$ such that $A + (-A) = A-A = 0$, where $0 \in \mathbb{R}$, rather than a fuzzy zero, which is not unique. 

Kosinski~\cite{Kosinski2003} defined an \textit{ordered fuzzy number} as follows:
\begin{defin}
An \textit{ordered fuzzy number} is an ordered pair $A = (\mu_{A}^{\uparrow}, \mu_{A}^{\downarrow})$ where $\mu_{A}^{\uparrow}, \mu_{A}^{\downarrow}:[0,1] \to \mathbb{R}$ are continuous functions. 
\end{defin}

With $\mathcal{C}[0,1]$ as the space of continuous functions from $[0,1]$ to $\mathbb{R}$, the above definition is equivalent to taking $A \in \mathcal{C}[0,1] \times \mathcal{C}[0,1]$.This pair of functions effectively defines the endpoints of the level sets of $A$ for a specified $\alpha \in [0,1]$.  The orientation of the ordered fuzzy number is given by the ordered pair, where one travels "up" $\mu_{A}^{\uparrow}$ from $\alpha = 0$ to $\alpha = 1$ and "down" $\mu_{A}^{\downarrow}$ from $\alpha = 1$ to $\alpha = 0$. The words "up" and "down" refer to the orientation, and not whether $\mu_{A}^{\uparrow}$ or $\mu_{A}^{\downarrow}$ are increasing or decreasing. 

\begin{defin}[Increasing/Decreasing OFN]
An OFN is \textit{increasing} \cite{MB2013} if the orientation is consistent with the direction of the axis, and decreasing if the orientation is opposite to the direction of the axis. 
\end{defin}

\begin{figure}[H]
\begin{tikzpicture}
\begin{axis}[scale=0.9, xmin=-1, xmax=1.5, ymin=-1, ymax=3, grid = both,
axis x line=middle,thick,
axis y line=middle,tick style={draw=none},
xlabel=$\alpha$, ylabel=$\mathbb{R}$,
every axis y label/.style={
    at={(ticklabel* cs:1.00)},
    anchor=west,
},
every axis x label/.style={
    at={(ticklabel* cs:1.00)},
    anchor=west,
}
]
\addplot[domain=0:1,samples=100,blue,
postaction={decorate, decoration={markings,
mark=at position 0.31 with {\arrow{>};},
mark=at position 0.71 with {\arrow{>};},
      }}
        ]{sqrt(x)}node[right,pos=0.9]{};
\node[below,text=blue,font=\large] at (0.6,0.7) {$\mu_{A}^{\uparrow} = \sqrt{\alpha}$};
\addplot[domain=0:1,samples=100,blue,
postaction={decorate, decoration={markings,
mark=at position 0.31 with {\arrow{<};},
mark=at position 0.71 with {\arrow{<};},
    }}
        ]{2-x}node[right,pos=0.9]{};
\node[below,text=blue,font=\large] at (0.6,2.5) {$\mu_{A}^{\downarrow} = 2-\alpha$};
\end{axis}
\end{tikzpicture}
\caption{An OFN $A = (\sqrt{\alpha}, 2-\alpha)$}
\label{fig: OFN 1}
\end{figure}
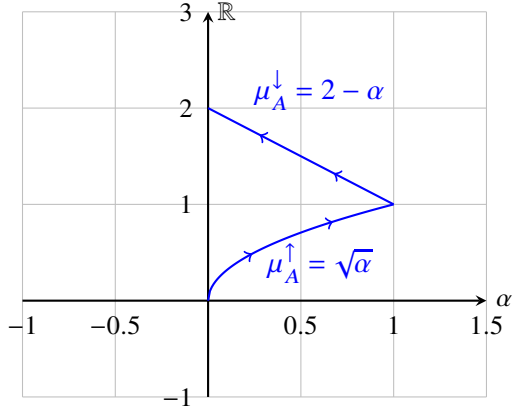

\begin{example}
 Let $\mu_{A}^{\uparrow}(\alpha) = \sqrt{\alpha}$, so that $\mu_{A}^{\uparrow}: [0,1] \to [0,1]$, and $\mu_{A}^{\downarrow}(\alpha) = 2-\alpha$, so that $\mu_{A}^{\downarrow}:[0,1] \to [1,2]$ as shown in Figure~\ref{fig: OFN 1}. Then $A = (\mu_{A}^{\uparrow}, \mu_{A}^{\downarrow})$ has level sets $A_{\alpha} = [\sqrt{\alpha}, 2-\alpha]$, giving $\text{supp}(A) = [0,2]$, $A_{1/4} = [1/2, 7/4]$, and $A_{1} = \{1\}$. 
 \end{example}

Kosinski never actually requires his OFNs to also be classical fuzzy numbers. It is possible to have an OFN $A = (\mu_{A}^{\uparrow}, \mu_{A}^{\downarrow})$ that does not fit the definition of a classical fuzzy number. Instead, he defines a certain class of \textit{proper} ordered fuzzy numbers as follows:

\begin{defin}[Proper OFN]
From ~\cite{Kacprzak2011}, an increasing OFN is \textit{proper} if $\mu_{A}^{\uparrow}$ is increasing, $\mu_{A}^{\downarrow}$ is decreasing, and $\mu_{A}^{\uparrow} \leq \mu_{A}^{\downarrow}$ for all $\alpha \in [0,1]$. A decreasing OFN is \textit{proper} if $\mu_{A}^{\downarrow}$ is increasing, $\mu_{A}^{\uparrow}$ is decreasing, and $\mu_{A}^{\downarrow} \leq \mu_{A}^{\uparrow}$ for all $\alpha \in [0,1]$.
\end{defin}

A proper OFN possesses a classical membership function. An OFN that is not proper is called \textit{improper}, and does not possess a classical membership function.

If an OFN $A = (\mu_{A}^{\uparrow}, \mu_{A}^{\downarrow})$ is proper, its membership function representation is given by 

\[\mu_{A}(x) = \begin{cases} (\mu_{A}^{\uparrow})^{-1}(x), & x \in [\mu_{A}^{\uparrow}(0), \mu_{A}^{\uparrow}(1)) \\ (\mu_{A}^{\downarrow})^{-1}(x), & x \in (\mu_{A}^{\downarrow}(1), \mu_{A}^{\downarrow}(0)] \\ 1, & x \in [\mu_{A}^{\uparrow}(1), \mu_{A}^{\downarrow}(1)]\end{cases}\]

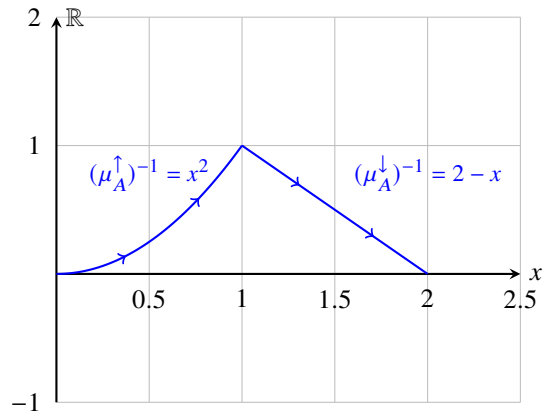
\begin{figure}[H]
\begin{tikzpicture}
\begin{axis}[scale=0.9, xmin=0, xmax=2.5, ymin=-1, ymax=2, grid = both,
axis x line=middle,thick,
axis y line=middle,tick style={draw=none},
xlabel=$x$, ylabel=$\mathbb{R}$,
every axis y label/.style={
    at={(ticklabel* cs:1.00)},
    anchor=west,
},
every axis x label/.style={
    at={(ticklabel* cs:1.00)},
    anchor=west,
}
]
\addplot[domain=0:1,samples=100,blue,
postaction={decorate, decoration={markings,
mark=at position 0.31 with {\arrow{>};},
mark=at position 0.71 with {\arrow{>};},
      }}
        ]{x^2}node[right,pos=0.9]{};
\node[below,text=blue,font=\small] at (0.5,1) {$(\mu_{A}^{\uparrow})^{-1} = x^{2}$};
\addplot[domain=1:2,samples=100,blue,
postaction={decorate, decoration={markings,
mark=at position 0.31 with {\arrow{>};},
mark=at position 0.71 with {\arrow{>};},
    }}
        ]{2-x}node[right,pos=0.9]{};
\node[below,text=blue,font=\small] at (2,1) {$(\mu_{A}^{\downarrow})^{-1} = 2-x$};
\end{axis}
\end{tikzpicture}
\caption{The membership function representation of the OFN $A = (\sqrt{\alpha}, 2-\alpha)$}
\label{fig: OFN 1 membership}
\end{figure}

For our example above, the membership function for $A = (\sqrt{\alpha}, 2-\alpha)$ is 

\[\mu_{A}(x) = \begin{cases} x^{2}, & x \in [0,1] \\
2-x, & x \in [1,2]
\end{cases}\]

If we swap the positions of the "up" and "down" functions, we create the same OFN shape with the opposite orientation. For example, let $B = (2-y,\sqrt{y})$. Then $B$ gives the same fuzzy set as $A$ above, but the orientation is reversed, so that we traverse the curve via the linear function first from $(2,0)$ to $(1,1)$, then move "down" $\sqrt{y}$ from $(1,1)$ to $(0,0)$. 

\begin{figure}[H]
\begin{tikzpicture}
\begin{axis}[scale=0.9, xmin=-1, xmax=1.5, ymin=-1, ymax=3, grid = both,
axis x line=middle,thick,
axis y line=middle,tick style={draw=none},
xlabel=$\alpha$, ylabel=$\mathbb{R}$,
every axis y label/.style={
    at={(ticklabel* cs:1.00)},
    anchor=west,
},
every axis x label/.style={
    at={(ticklabel* cs:1.00)},
    anchor=west,
}
]
\addplot[domain=0:1,samples=100,blue,
postaction={decorate, decoration={markings,
mark=at position 0.31 with {\arrow{<};},
mark=at position 0.71 with {\arrow{<};},
      }}
        ]{sqrt(x)}node[right,pos=0.9]{};
\node[below,text=blue,font=\large] at (0.6,0.7) {$\mu_{A}^{\downarrow} = \sqrt{\alpha}$};
\addplot[domain=0:1,samples=100,blue,
postaction={decorate, decoration={markings,
mark=at position 0.31 with {\arrow{>};},
mark=at position 0.71 with {\arrow{>};},
    }}
        ]{2-x}node[right,pos=0.9]{};
\node[below,text=blue,font=\large] at (0.6,2.5) {$\mu_{A}^{\uparrow} = 2-\alpha$};
\end{axis}
\end{tikzpicture}
\caption{The OFN $B = (2-\alpha,\sqrt{\alpha})$}
\label{fig: OFN 2}
\end{figure}
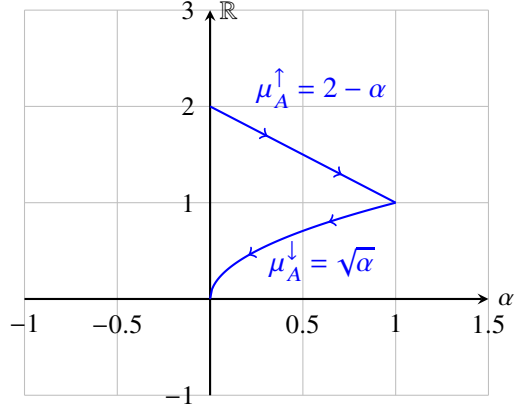

In ~\cite{Kosinski2003b}, Kosinski defines four arithmetic operations ($+,-,\cdot, /$) as pairwise operations on the components of the OFN. That is, for any above operation $\star$, 
\[C = A \star B = (\mu_{A}^{\uparrow}\star \mu_{B}^{\uparrow}, \mu_{A}^{\downarrow}\star \mu_{B}^{\downarrow})\]

With the added the property of orientation to fuzzy numbers, these ordered fuzzy numbers were intended to produce a unique additive inverse for each OFN. 

The additive inverse of an OFN $A = (\mu_{A}^{\uparrow}, \mu_{A}^{\downarrow})$ is $-A = (-\mu_{A}^{\uparrow},-\mu_{A}^{\downarrow})$. Thus, according to Kosinski's arithmetic operation, 
\[A + (-A) = (0,0)\]
which is a crisp $0$. One will also notice that $-A = -1\cdot A = -1(\mu_{A}^{\uparrow}, \mu_{A}^{\downarrow})$, where $-1 \in \mathbb{R}$. ~\cite{Kosinski2003b} shows then that the set $\mathbb{K}$ of all ordered fuzzy numbers is a vector space over $\mathbb{R}$ under this defined addition. 

\subsection{Drawbacks and Improper OFNs}
\label{subsec: improper OFNs}

The set $\mathbb{K}$ of all ordered fuzzy numbers only requires continuity of $\mu_{A}^{\uparrow}, \mu_{A}^{\downarrow}$. With additional requirements as stated above, an OFN may be a proper OFN, meaning $\mu_{A}^{\uparrow}, \mu_{A}^{\downarrow}$ can be inverted and the OFN $A$ has a membership function that fits the definition of a classical fuzzy number. 

Any OFN that is not proper is improper. From the definition above of a proper OFN, we may deduce that an improper OFN has at least one of the following properties.

\begin{enumerate}[(i)]
\item At least one of $\mu_{A}^{\uparrow}, \mu_{A}^{\downarrow}$ is not monotonic,
\item $\mu_{A}^{\uparrow}, \mu_{A}^{\downarrow}$ are either both increasing or both decreasing,
\item $|\text{range}(\mu_{A}^{\uparrow}) \cap \text{range}(\mu_{A}^{\downarrow})| > 1$, where $| \cdot |$ denotes the cardinality of a set.
\end{enumerate}

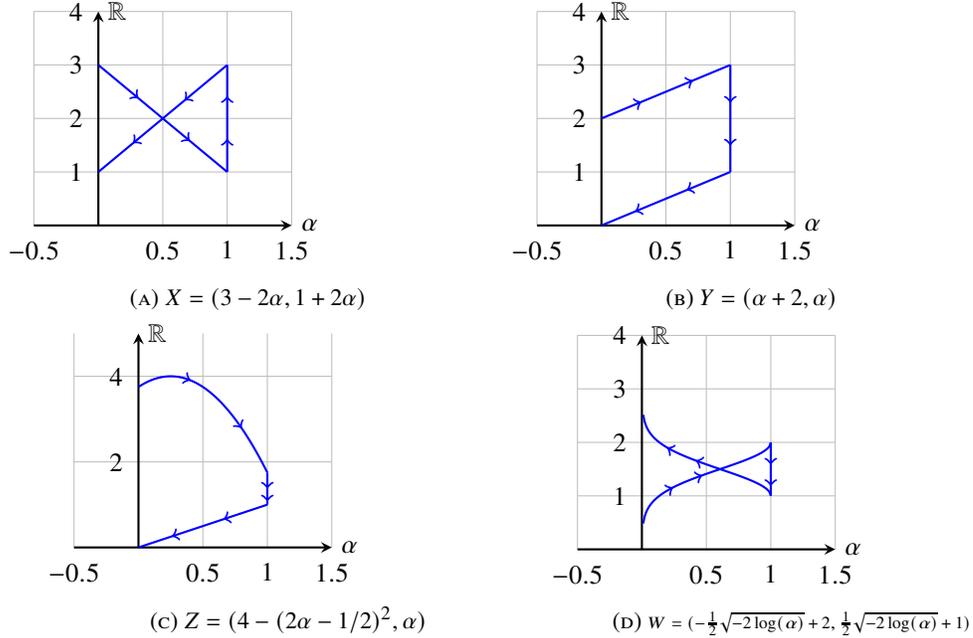
\begin{figure}[H]
\begin{subfigure}{0.4\linewidth}
\begin{tikzpicture}
\begin{axis}[scale=0.5, xmin=-0.5, xmax=1.5, ymin=0, ymax=4, grid = both,
axis x line=middle,thick,
axis y line=middle,tick style={draw=none},
xlabel=$\alpha$, ylabel=$\mathbb{R}$,
every axis y label/.style={
    at={(ticklabel* cs:1.00)},
    anchor=west,
},
every axis x label/.style={
    at={(ticklabel* cs:1.00)},
    anchor=west,
}
]
\addplot[domain=0:1,samples=100,blue,
postaction={decorate, decoration={markings,
mark=at position 0.31 with {\arrow{<};},
mark=at position 0.71 with {\arrow{<};},
      }}
        ]{2*x + 1}node[right,pos=0.9]{};
\node[below,text=blue,font=\large] at (0.6,0.7) {};
\addplot[domain=0:1,samples=100,blue,
postaction={decorate, decoration={markings,
mark=at position 0.31 with {\arrow{>};},
mark=at position 0.71 with {\arrow{>};},
    }}
        ]{3-2*x}node[right,pos=0.9]{};
\node[below,text=blue,font=\large] at (0.6,2.5) {};
\addplot[domain=0:1,samples=100,blue,
postaction={decorate, decoration={markings,
mark=at position 0.31 with {\arrow{>};},
mark=at position 0.71 with {\arrow{>};},
    }}
        ] coordinates{(1,1)(1,3)};
\end{axis}
\end{tikzpicture}
\caption{$X = (3-2\alpha, 1+2\alpha)$}
\label{subfig: improper 1}
\end{subfigure}
\begin{subfigure}{0.4\linewidth}
\begin{tikzpicture}
\begin{axis}[scale=0.5, xmin=-0.5, xmax=1.5, ymin=0, ymax=4, grid = both,
axis x line=middle,thick,
axis y line=middle,tick style={draw=none},
xlabel=$\alpha$, ylabel=$\mathbb{R}$,
every axis y label/.style={
    at={(ticklabel* cs:1.00)},
    anchor=west,
},
every axis x label/.style={
    at={(ticklabel* cs:1.00)},
    anchor=west,
}
]
\addplot[domain=0:1,samples=100,blue,
postaction={decorate, decoration={markings,
mark=at position 0.31 with {\arrow{<};},
mark=at position 0.71 with {\arrow{<};},
      }}
        ]{x}node[right,pos=0.9]{};
\node[below,text=blue,font=\large] at (0.6,0.7) {};
\addplot[domain=0:1,samples=100,blue,
postaction={decorate, decoration={markings,
mark=at position 0.31 with {\arrow{>};},
mark=at position 0.71 with {\arrow{>};},
    }}
        ]{2+x}node[right,pos=0.9]{};
\node[below,text=blue,font=\large] at (0.6,2.5) {};
\addplot[domain=0:1,samples=100,blue,
postaction={decorate, decoration={markings,
mark=at position 0.31 with {\arrow{<};},
mark=at position 0.71 with {\arrow{<};},
    }}
        ] coordinates{(1,1)(1,3)};
\end{axis}
\end{tikzpicture}
\caption{$Y = (\alpha+2, \alpha)$}
\label{subfig: improper 2}
\end{subfigure}
\newline
\begin{subfigure}{0.4\linewidth}
\begin{tikzpicture}
\begin{axis}[scale=0.5, xmin=-0.5, xmax=1.5, ymin=0, ymax=5, grid = both,
axis x line=middle,thick,
axis y line=middle,tick style={draw=none},
xlabel=$\alpha$, ylabel=$\mathbb{R}$,
every axis y label/.style={
    at={(ticklabel* cs:1.00)},
    anchor=west,
},
every axis x label/.style={
    at={(ticklabel* cs:1.00)},
    anchor=west,
}
]
\addplot[domain=0:1,samples=100,blue,
postaction={decorate, decoration={markings,
mark=at position 0.31 with {\arrow{<};},
mark=at position 0.71 with {\arrow{<};},
      }}
        ]{x}node[right,pos=0.9]{};
\node[below,text=blue,font=\large] at (0.6,0.7) {};
\addplot[domain=0:1,samples=100,blue,
postaction={decorate, decoration={markings,
mark=at position 0.31 with {\arrow{>};},
mark=at position 0.71 with {\arrow{>};},
    }}
        ]{-(2*x-0.5)^2 + 4}node[right,pos=0.9]{};
\node[below,text=blue,font=\large] at (0.6,2.5) {};
\addplot[domain=0:1,samples=100,blue,
postaction={decorate, decoration={markings,
mark=at position 0.31 with {\arrow{<};},
mark=at position 0.71 with {\arrow{<};},
    }}
        ] coordinates{(1,1)(1,7/4)};
\end{axis}
\end{tikzpicture}
\caption{$Z = (4-(2\alpha-1/2)^2, \alpha)$}
\label{subfig: improper 3}
\end{subfigure}
\begin{subfigure}{0.4\linewidth}
\begin{tikzpicture}
\begin{axis}[scale=0.5, xmin=-0.5, xmax=1.5, ymin=0, ymax=4, grid = both,
axis x line=middle,thick,
axis y line=middle,tick style={draw=none},
xlabel=$\alpha$, ylabel=$\mathbb{R}$,
every axis y label/.style={
    at={(ticklabel* cs:1.00)},
    anchor=west,
},
every axis x label/.style={
    at={(ticklabel* cs:1.00)},
    anchor=west,
}
]
\addplot[domain=0:1,samples=100,blue,
postaction={decorate, decoration={markings,
mark=at position 0.31 with {\arrow{<};},
mark=at position 0.51 with {\arrow{<};},
      }}
        ]{0.5*sqrt(-2*ln(x))+1}node[right,pos=0.9]{};
\node[below,text=blue,font=\large] at (0.6,0.7) {};
\addplot[domain=0:1,samples=100,blue,
postaction={decorate, decoration={markings,
mark=at position 0.31 with {\arrow{>};},
mark=at position 0.51 with {\arrow{>};},
    }}
        ]{-0.5*sqrt(-2*ln(x))+2}node[right,pos=0.9]{};
\node[below,text=blue,font=\large] at (0.6,2.5) {};
\addplot[domain=0:1,samples=100,blue,
postaction={decorate, decoration={markings,
mark=at position 0.31 with {\arrow{<};},
mark=at position 0.71 with {\arrow{<};},
    }}
        ] coordinates{(1,1)(1,2)};
\end{axis}
\end{tikzpicture}
\caption{\tiny$W = (-\tfrac{1}{2}\sqrt{-2\log(\alpha)} + 2, \tfrac{1}{2}\sqrt{-2\log(\alpha)}+1)$}
\label{subfig: improper 4}
\end{subfigure}
\caption{Selected improper OFNs}
\label{fig: improper OFN}
\end{figure}

Figure~\ref{fig: improper OFN} shows selected examples of improper OFNs. Figure~\ref{subfig: improper 1} is an example of (iii), Figure~\ref{subfig: improper 2} results from (ii), and Figure~\ref{subfig: improper 3} results from condition (i), though it should be noted that all of these pathological conditions may be present in a single OFN. Improper OFNs lack a membership function, which has drawbacks in practice. An improper OFN has at least one point $x$ with more than one membership level, which is nonsensical when the OFN is intended for classical use as a fuzzy set. 

A particular drawback to Kosinski's arithmetic operations on OFNs is that the sum or product of two proper OFNs may yield an improper OFN, and the impropriety is difficult to control, especially regarding the pathology of non-monotonicity in the up and down functions. 

For an example consider the two proper OFNs $A = (x, 2-x)$, $B = (1-2x, -3)$. Then
\[A+B = (1-x, -1-x)\]
This is an improper fuzzy number, since we have that $(A+B)_{1/2} = [-3/2, 1/2]$ and yet $(A+B)_{1} = [-2,0]$. $(A+B)_{1} \not\subset (A+B)_{1/2}$, and thus $(A+B)$ is an improper fuzzy number, despite being in $\mathbb{K}$. 
Kosinski's operations present a challenge regarding improper OFNs and the control of such issues over arithmetic operations, particularly repeated operations. If $\hat{\mathbb{K}}$ denotes the set of proper OFNs, then none of Kosinski's operations will close $\hat{\mathbb{K}}$, and moreover, pathologies become difficult to control with repeated operations on a set as broadly defined as the set of all OFNs $\mathbb{K}$.~\cite{MB2013} mentions this issue in repeated applications of Kosinski's arithmetic on OFNs when using them to estimate regression coefficients for fuzzy autoregressive models. 

Nonetheless, it is desired to have as strong an algebraic structure as possible when defining operations over OFNs, while attempting to maintain applicability to a multitude of situations. Piasecki~\cite{Piasecki2018} attempted to define a different collection of arithmetic operations that would close the set of proper OFNs $\hat{\mathbb{K}}$ but unfortunately failed. These operations are defined for trapezoidal OFNs thus:

 Let $X = (\mu_{X}^{\uparrow}, \mu_{X}^{\downarrow})$, $Y = (\mu_{Y}^{\uparrow}, \mu_{Y}^{\downarrow})$. Let $a_{\cdot} = \mu_{\cdot}^{\uparrow}(0)$, $b_{\cdot} = \mu_{\cdot}^{\uparrow}(1)$, $c_{\cdot} = \mu_{\cdot}^{\downarrow}(1)$, $d_{\cdot} = \mu_{\cdot}^{\downarrow}(0)$, where $\cdot = X,Y$. Then any arithmetic operation $\star$ on $\mathbb{R}$ is extended to the set $\hat{\mathbb{K}}$ of proper OFNs as follows: 

For OFN $W = X \tilde{\star} Y$,
\begin{equation}
\begin{aligned}
a_{W} &= \begin{cases} \min\{a_{X}\star a_{Y}, b_{X} \star b_{Y}\}, & (b_{X} \star b_{Y} < c_{X} \star c_{Y}) \text{ or } (b_{X}\star b_{Y} = c_{X} \star c_{Y} \text{ and } a_{X} \star a_{Y} \leq d_{X} \star d_{Y}) \\
\max\{a_{X}\star a_{Y}, b_{X} \star b_{Y}\}, & (b_{X} \star b_{Y} > c_{X} \star c_{Y}) \text{ or } (b_{X}\star b_{Y} = c_{X} \star c_{Y} \text{ and } a_{X} \star a_{Y} > d_{X} \star d_{Y})\end{cases} \\
b_{W} &= b_{X} \star b_{Y} \\
c_{W} &= c_{X} \star c_{Y} \\
d_{W} &= \begin{cases} \max\{d_{X}\star d_{Y}, c_{X} \star c_{Y}\}, & (b_{X} \star b_{Y} < c_{X} \star c_{Y}) \text{ or } (b_{X}\star b_{Y} = c_{X} \star c_{Y} \text{ and } a_{X} \star a_{Y} \leq d_{X} \star d_{Y}) \\
\min\{d_{X}\star d_{Y}, c_{X} \star c_{Y}\}, & (b_{X} \star b_{Y} > c_{X} \star c_{Y}) \text{ or } (b_{X}\star b_{Y} = c_{X} \star c_{Y} \text{ and } a_{X} \star a_{Y} > d_{X} \star d_{Y})\end{cases} \\
\mu_{W}^{\uparrow}(\alpha) &= \begin{cases} \mu_{X}^{\uparrow}(\alpha) \star \mu_{Y}^{\uparrow}(\alpha), & a_{W} \neq b_{W} \\
b_{W}, & a_{W} = b_{W} \end{cases} \\
\mu_{W}^{\downarrow}(\alpha) &= \begin{cases} \mu_{X}^{\downarrow}(\alpha) \star \mu_{Y}^{\downarrow}(\alpha), & c_{W} \neq d_{W} \\
c_{W}, & c_{W} = d_{W} 
\end{cases}
\end{aligned}
\end{equation}

The claim that these revised operations close the space of $\hat{\mathbb{K}}$ is incorrect. We provide the following counterexample to  Theorem 5 of ~\cite{Piasecki2018}:

Let $X = (-1,-x)$, $Y = (2,4x-2)$. Then $X,Y \in \hat{\mathbb{K}}$. We look at the operation "multiplication" under Piasecki's definition. Then letting $A = X \mathrel{\tilde{\cdot}} Y$, 

\begin{align*}
a_{A} &= \min\{a_{X}\cdot a_{Y}, b_{X} \cdot b_{Y}\} = -2 \\
b_{A} &= b_{X} \cdot b_{Y} = -2 \\
c_{A} &= c_{X} \cdot c_{Y} = -2 \\
d_{A} &= \max\{d_{X}\star d_{Y}, c_{X} \star c_{Y}\} = 0
\end{align*}
Since $a_{A} = b_{A}$, $\mu_{A}^{\uparrow}(\alpha) = b_{A} = -2$. Since $c_{A} \neq d_{A}$, $\mu_{A}^{\downarrow}(\alpha) = \mu_{X}^{\downarrow}(\alpha)\cdot \mu_{Y}^{\downarrow}(\alpha) = -x(4x-2)$. However, notice that $\mu_{A}^{\downarrow}$ is not a monotonic function. In particular, $\mu_{A}^{\downarrow}(0) = 0$, and $\mu_{A}^{\downarrow}(1/2) = 0$, which implies that the point 0 has membership levels 0 and 1/2 simultaneously. Thus, $\mu_{A}^{\downarrow}$ cannot be inverted to yield a proper membership function for a fuzzy number, and Piasecki's modification of multiplication does not close $\hat{\mathbb{K}}$. 

We now present the concept of typed OFNs and define arithmetic operations over smaller subsets of OFNs such that we construct rings of typed OFNs. Our typed OFNs present a more restricted view of OFNs than Kosinski~\cite{Kosinski2003}, but yield an advantage in controlling pathologies resulting from improper OFNs, and remain closed under a particular type, correcting the issues of~\cite{Piasecki2018}. We detail this in the next section. 

\section{Typed OFNs}

In order to define typed OFNs, we present some preliminary definitions. It should be assumed that we are always working with continuous, real-valued functions.

\begin{defin}[Elementary function]
An \textit{elementary function}~\cite{Spivakcalc} is one which can be obtained by addition, multiplication, division, and composition from the rational functions, the trigonometric functions and their inverses, and the functions $\log$ and $\exp$. 
\end{defin}
\begin{defin}[Parent Function]
A \textit{parent function} is a single-term elementary function.
\end{defin}

\begin{defin}[Base function]
An OFN base function $h(x)$ is either a parent function or combination of parent functions obtained through a combination of inversion and/or composition. 
\end{defin}

In other words, base functions may only be parent functions, compositions of parent functions, inverses of parent functions, or combinations of composition and inversion of parent functions. 

Some examples of parent functions include
\begin{enumerate}[(a)]
\item For the class of linear functions with the standard form $g(x) = mx+b$, both $x$ and $-x$ are parent functions.
\item For the class of exponential functions, $e^{x}$ and $-e^{x}$ are parent functions.
\item For the class of $n$th degree polynomials, $x^{n}$ and $-x^{n}$ are parent functions for $n \in \mathbb{N}$
\item For the class of natural logarithms, $\log(x)$ is a parent function.
\end{enumerate}

Some examples of base functions include

\begin{enumerate}[(a)]
\item $h(x) = e^{-x}$, given as the composition of $e^{x}$ with $-x$.
\item $h(x) = \sqrt{-2\log(x)}$, given as the inverse of $e^{-x^{2}}$, which itself is built from the composition of $e^{x}$ with $-x^{2}$. 
\end{enumerate}

We note that constant functions should be treated separately, and define them thusly:
\begin{defin}[Typeless function]
A constant function will be called \textit{typeless}.
\end{defin} 

We are now ready to define a typed OFN. 

\begin{defin}[True-Typed OFN]
Let $h(x)$ be a fixed base function. We say an OFN $A = (\mu_{A}^{\uparrow}, \mu_{A}^{\downarrow})$ is an \textit{true OFN of type h over $\mathbb{R}$} or an $h-$OFN if both $\mu_{A}^{\uparrow}, \mu_{A}^{\downarrow}$ have the form $f(x) = a h(x) + b$
for $a, b \in \mathbb{R}$, where $f: [0,1] \to \mathbb{R}$.
\end{defin}

We now give an equivalent representation of a typed OFN upon which we will define our operations.

\begin{defin}[Typed OFN]
For a fixed base function $h$, we say an OFN $A = (\mu_{A}^{\uparrow}, \mu_{A}^{\downarrow})$ is a \textit{typed OFN} if at least one of $\mu_{A}^{\uparrow}, \mu_{A}^{\downarrow}$ has the form $f(x) = a h(x) + b$, $a,b,\in \mathbb{R}$, with the other possibly typeless.
\end{defin}

\begin{defin}[Essential tuple]
For a given $h-$OFN $A = (\mu_{A}^{\uparrow}, \mu_{A}^{\downarrow}) = (a^{\uparrow} h(x) + b^{\uparrow}, a^{\downarrow} h(x) + b^{\downarrow})$, the tuple $(a^{\uparrow}, b^{\uparrow}, a^{\downarrow}, b^{\downarrow})$, with $a^{\uparrow}, b^{\uparrow}, a^{\downarrow}, b^{\downarrow} \in \mathbb{R}$ is called the \textit{essential tuple}.
\end{defin}

\begin{defin}[Rectangular or Typeless OFN]
An OFN $A = (\mu_{A}^{\uparrow}, \mu_{A}^{\downarrow})$ where both of $\mu_{A}^{\uparrow}, \mu_{A}^{\downarrow}$ are typeless is called a \textit{rectangular} or \textit{typeless} OFN, and has essential tuple $(0,b^{\uparrow}, 0, b^{\downarrow})$, for $b^{\uparrow}, b^{\downarrow} \in \mathbb{R}$.
\end{defin}

\begin{defin}[OFN of mixed type]
If $\mu_{A}^{\uparrow}$ and $\mu_{A}^{\downarrow}$ are not constructed from the same base function, then we call $A$ an \textit{OFN of mixed type}.
\end{defin}

\textit{Remark: We do not define essential tuples for mixed-type OFNs.}


\begin{figure}[H]
\begin{subfigure}{0.4\linewidth}
\begin{tikzpicture}
\begin{axis}[scale=0.5, xmin=-0.5, xmax=1.5, ymin=0, ymax=4, grid = both,
axis x line=middle,thick,
axis y line=middle,tick style={draw=none},
xlabel=$\alpha$, ylabel=$\mathbb{R}$,
every axis y label/.style={
    at={(ticklabel* cs:1.00)},
    anchor=west,
},
every axis x label/.style={
    at={(ticklabel* cs:1.00)},
    anchor=west,
}
]
\addplot[domain=0:1,samples=100,blue,
postaction={decorate, decoration={markings,
mark=at position 0.31 with {\arrow{<};},
mark=at position 0.71 with {\arrow{<};},
      }}
        ]{x}node[right,pos=0.9]{};
\node[below,text=blue,font=\large] at (0.6,0.7) {};
\addplot[domain=0:1,samples=100,blue,
postaction={decorate, decoration={markings,
mark=at position 0.31 with {\arrow{>};},
mark=at position 0.71 with {\arrow{>};},
    }}
        ]{3-x}node[right,pos=0.9]{};
\node[below,text=blue,font=\large] at (0.6,2.5) {};
\addplot[domain=0:1,samples=100,blue,
postaction={decorate, decoration={markings,
mark=at position 0.31 with {\arrow{<};},
mark=at position 0.71 with {\arrow{<};},
    }}
        ] coordinates{(1,1)(1,2)};
\end{axis}
\end{tikzpicture}
\caption{Truly trapezoidal OFN\\ $A = (3-x,x)$}
\label{fig: typed OFN 1}
\end{subfigure}
\begin{subfigure}{0.4\linewidth}
\begin{tikzpicture}
\begin{axis}[scale=0.5, xmin=-0.5, xmax=1.5, ymin=0, ymax=4, grid = both,
axis x line=middle,thick,
axis y line=middle,tick style={draw=none},
xlabel=$\alpha$, ylabel=$\mathbb{R}$,
every axis y label/.style={
    at={(ticklabel* cs:1.00)},
    anchor=west,
},
every axis x label/.style={
    at={(ticklabel* cs:1.00)},
    anchor=west,
}
]
\addplot[domain=0:1,samples=100,blue,
postaction={decorate, decoration={markings,
mark=at position 0.31 with {\arrow{<};},
mark=at position 0.71 with {\arrow{<};},
      }}
        ]{sqrt(-2*ln(x)) + 1}node[right,pos=0.9]{};
\node[below,text=blue,font=\large] at (0.6,0.7) {};
\addplot[domain=0:1,samples=100,blue,
postaction={decorate, decoration={markings,
mark=at position 0.31 with {\arrow{>};},
mark=at position 0.71 with {\arrow{>};},
    }}
        ]{1}node[right,pos=0.9]{};
\node[below,text=blue,font=\large] at (0.6,2.5) {};
\end{axis}
\end{tikzpicture}
\caption{Gaussian OFN, with one side rectangular \\ $B = (1, \sqrt{-2\log(x)}+1)$}
\label{fig: typed OFN 2}
\end{subfigure}
\newline
\begin{subfigure}{0.4\linewidth}
\begin{tikzpicture}
\begin{axis}[scale=0.5, xmin=-0.5, xmax=1.5, ymin=0, ymax=5, grid = both,
axis x line=middle,thick,
axis y line=middle,tick style={draw=none},
xlabel=$\alpha$, ylabel=$\mathbb{R}$,
every axis y label/.style={
    at={(ticklabel* cs:1.00)},
    anchor=west,
},
every axis x label/.style={
    at={(ticklabel* cs:1.00)},
    anchor=west,
}
]
\addplot[domain=0:1,samples=100,blue,
postaction={decorate, decoration={markings,
mark=at position 0.31 with {\arrow{>};},
mark=at position 0.71 with {\arrow{>};},
      }}
        ]{2}node[right,pos=0.9]{};
\node[below,text=blue,font=\large] at (0.6,0.7) {};
\addplot[domain=0:1,samples=100,blue,
postaction={decorate, decoration={markings,
mark=at position 0.31 with {\arrow{<};},
mark=at position 0.71 with {\arrow{<};},
    }}
        ]{1}node[right,pos=0.9]{};
\node[below,text=blue,font=\large] at (0.6,2.5) {};
\addplot[domain=0:1,samples=100,blue,
postaction={decorate, decoration={markings,
mark=at position 0.31 with {\arrow{<};},
mark=at position 0.71 with {\arrow{<};},
    }}
        ] coordinates{(1,1)(1,2)};
\end{axis}
\end{tikzpicture}
\caption{Rectangular (typeless) OFN \\$C = (2,1)$}
\label{fig: typed OFN 3}
\end{subfigure}
\begin{subfigure}{0.4\linewidth}
\begin{tikzpicture}
\begin{axis}[scale=0.5, xmin=-0.5, xmax=1.5, ymin=0, ymax=4, grid = both,
axis x line=middle,thick,
axis y line=middle,tick style={draw=none},
xlabel=$\alpha$, ylabel=$\mathbb{R}$,
every axis y label/.style={
    at={(ticklabel* cs:1.00)},
    anchor=west,
},
every axis x label/.style={
    at={(ticklabel* cs:1.00)},
    anchor=west,
}
]
\addplot[domain=0:1,samples=100,blue,
postaction={decorate, decoration={markings,
mark=at position 0.31 with {\arrow{<};},
mark=at position 0.51 with {\arrow{<};},
      }}
        ]{-x+3}node[right,pos=0.9]{};
\node[below,text=blue,font=\large] at (0.6,0.7) {};
\addplot[domain=0:1,samples=100,blue,
postaction={decorate, decoration={markings,
mark=at position 0.31 with {\arrow{>};},
mark=at position 0.51 with {\arrow{>};},
    }}
        ]{0.5*ln(x)+2}node[right,pos=0.9]{};
\node[below,text=blue,font=\large] at (0.6,2.5) {};
\end{axis}
\end{tikzpicture}
\caption{OFN of mixed type \\ $D = (\tfrac{1}{2}\log(x) + 2, 3-x)$}
\label{fig: typed OFN 4}
\end{subfigure}
\caption{Examples of Typed OFNs}
\label{fig: typed OFN}
\end{figure}
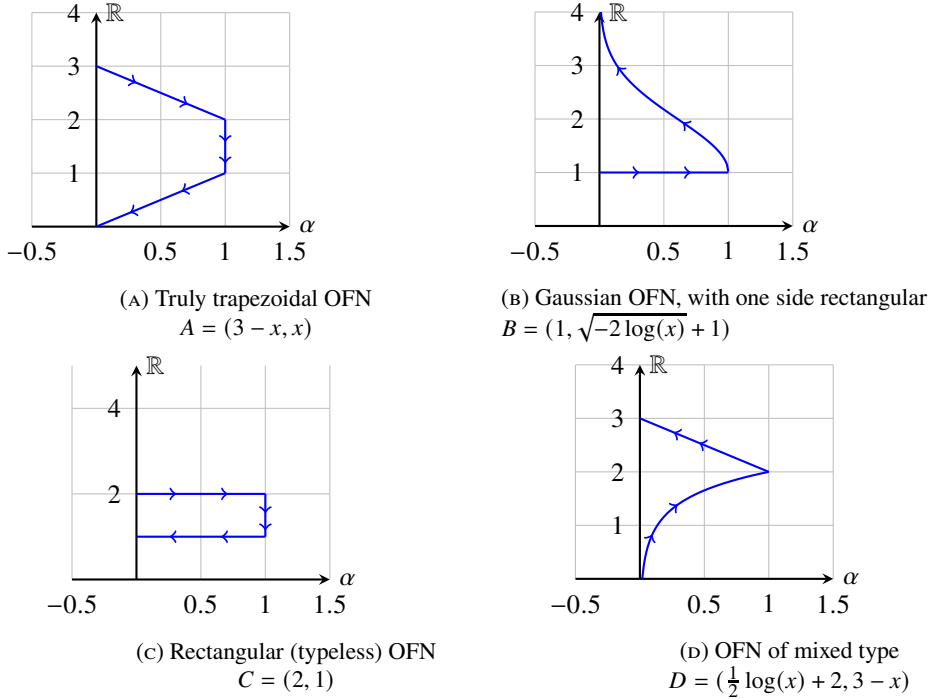

Figure~\ref{fig: typed OFN 1} is a true linear (or true trapezoidal) OFN with essential tuple $(-1,3,1,0)$; Figure~\ref{fig: typed OFN 2} represents a Gaussian OFN with $\mu_{A}^{\uparrow}$ typeless, with essential tuple $(0,1,1,1)$; Figure~\ref{fig: typed OFN 3} is a rectangular OFN with essential tuple $(0,2,0,1)$; and Figure~\ref{fig: typed OFN 4} is an OFN of mixed type and therefore does not have an essential tuple. 

We choose to restrict our attention to typed OFNs built from base functions which are monotonic on $\mathbb{R}$, and typeless base functions. This will allow us to control pathologies that may result from arithmetic on such OFNs. The next section defines arithmetic operations on the essential tuples of a set of typed OFNs to construct a commutative ring.

\section{Constructing rings of Typed OFNs}
 
We now present the general construction of rings of typed OFNs for a specified type. Let $\mathbb{K}_{h}$ denote the set of all $h-$OFNs for a given base function $h(x)$ (to include possibly one side of the OFN typeless), let $\mathbb{K}_{r}$ denote the set of rectangular or typeless OFNs. Additionally, denote $\mathbb{K}_{h,r} = \mathbb{K}_{h} \cup \mathbb{K}_{r}$.

We define arithmetic operations over the set of essential tuples of $\mathbb{K}_{h,r}$ in the following manner:

Let $\xi_{X} = (a_{X}^{\uparrow}, b_{X}^{\uparrow}, a_{X}^{\downarrow}, b_{X}^{\downarrow})$ be the essential tuple for $X \in \mathbb{K}_{h,r}$. For $X,Y \in \mathbb{K}_{h,r}$ and arithmetic operation $\star \in \left\{+,-,\cdot, / \right\}$ on $\mathbb{R}$, we define $A = X \inbox{$\star$} Y \in \mathbb{K}_{h,r}$ as the $h-$OFN with essential tuple 
\[\xi_{A} = (a_{X}^{\uparrow}\star a_{Y}^{\uparrow}, b_{X}^{\uparrow}\star b_{Y}^{\uparrow}, a_{X}^{\downarrow}\star a_{Y}^{\downarrow}, b_{X}^{\downarrow} \star b_{Y}^{\downarrow})\]

Thus, arithmetic operations take place component-wise over the essential tuples of a given set of typed OFNs.  Note that $\boxplus$
and $\boxdot$ are equivalent to the ordinary vector sum and Hadamard
product respectively on $\mathbb{R}^4$.
Arithmetic operations are only permitted between OFNs in the same typed set, and arithmetic operations on OFNs of mixed type or between a typed OFN and one of mixed type are not permitted or defined, as the essential tuple only uniquely defined an OFN after the type has been specified. 

\begin{theorem}
Let $h$ be a specified base function, and $\inbox{$\star$}$ the above defined operations on $\mathbb{K}_{h,r}$. Then $(\mathbb{K}_{h,r}, \boxplus, \boxdot)$ is a commutative ring.
\end{theorem}
\label{thm: ring}

\begin{proof}
An OFN $A =(\mu_{A}^{\uparrow}, \mu_{A}^{\downarrow}) \in \mathbb{K}_{h,r}$ is uniquely defined by its essential tuple $\xi_{A} \in \mathbb{R}$. for a specified $h$. Let $X, Y \in \mathbb{K}_{h,r}$ with essential tuples $\xi_{X}, \xi_{Y}$. Since $X \inbox{$\star$} Y$ produces an essential tuple in $\mathbb{K}_{h,r}$, the operations close $\mathbb{K}_{h,r}$.
We verify the ring axioms on the essential tuples. 

Let $X,Y,Z \in \mathbb{K}_{h,r}$ with essential tuples $\xi_{X} = (a_{X}^{\uparrow}, b_{X}^{\uparrow}, a_{X}^{\downarrow}, b_{X}^{\downarrow})$, $\xi_{Y} = (a_{Y}^{\uparrow}, b_{Y}^{\uparrow}, a_{Y}^{\downarrow}, b_{Y}^{\downarrow})$, $\xi_{Z} = (a_{Z}^{\uparrow}, b_{Z}^{\uparrow}, a_{Z}^{\downarrow}, b_{Z}^{\downarrow})$. 

\begin{enumerate}[(1)]
\item \textbf{Associativity of Addition:} 
Let $W \in \mathbb{K}_{h,r}$ with essential tuple $\xi_{W}$ such that $W = (X \boxplus Y) \boxplus Z$, and $U \in \mathbb{K}_{h,r}$ with essential tuple $\xi_{U}$ such that $U = X \boxplus (Y \boxplus Z)$. 
Then
\begin{align*}
\xi_{W} &= ((a_{X}^{\uparrow} + a_{Y}^{\uparrow}) + a_{Z}^{\uparrow}, (b_{X}^{\uparrow} + b_{Y}^{\uparrow}) + b_{Z}^{\uparrow}, (a_{X}^{\downarrow} + a_{Y}^{\downarrow}) + a_{Z}^{\downarrow}, (b_{X}^{\downarrow} + b_{Y}^{\downarrow}) + b_{Z}^{\downarrow}) \\
 &= (a_{X}^{\uparrow} + (a_{Y}^{\uparrow} + a_{Z}^{\uparrow}), b_{X}^{\uparrow} + (b_{Y}^{\uparrow} + b_{Z}^{\uparrow}), a_{X}^{\downarrow} + (a_{Y}^{\downarrow} + a_{Z}^{\downarrow}), b_{X}^{\downarrow} + (b_{Y}^{\downarrow} + b_{Z}^{\downarrow})) \\
 &= \xi_{U}
 \end{align*}
As the essential tuples uniquely determine an OFN, $W = U$.
 \item \textbf{Commutativity of Addition}: We have that 
 \begin{align*}
 \xi_{X} \boxplus \xi_{Y} &= (a_{X}^{\uparrow} + a_{Y}^{\uparrow}, b_{X}^{\uparrow} + b_{Y}^{\uparrow}, a_{X}^{\downarrow} + a_{Y}^{\downarrow}, b_{X}^{\downarrow} + b_{Y}^{\downarrow}) \\
 &=  (a_{Y}^{\uparrow} + a_{X}^{\uparrow}, b_{Y}^{\uparrow} + b_{X}^{\uparrow}, a_{Y}^{\downarrow} + a_{X}^{\downarrow}, b_{Y}^{\downarrow} + b_{X}^{\downarrow}) \\
 &= \xi_{Y} \boxplus \xi_{X}
 \end{align*} 
 where the second equality comes from commutativity of $+$ on $\mathbb{R}$.
 \item \textbf{Additive Identity:} The crisp $0 \in \mathbb{R}$ whose essential tuple is $(0,0,0,0)$ is the additive identity for $\mathbb{K}_{h,r}$, since addition is performed componentwise over $\mathbb{R}$, and $0$ is the additive identity of $\mathbb{R}$. 
 \item \textbf{Additive Inverse:} We have that $-X$ with $\xi_{-X} = (-a_{X}^{\uparrow}, -b_{X}^{\uparrow}, -a_{X}^{\downarrow}, -b_{X}^{\downarrow})$ is the unique additive inverse of $X$, as $X \boxplus -X$ yields the essential tuple $(0,0,0,0)$, which is the additive identity.
 \item \textbf{Associativity of Multiplication:}
 Let $W \in \mathbb{K}_{h,r}$ with essential tuple $\xi_{W}$ such that $W = (X \boxtimes Y) \boxtimes Z$, and $U \in \mathbb{K}_{h,r}$ with essential tuple $\xi_{U}$ such that $U = X \boxtimes (Y \boxtimes Z)$. 
Then
\begin{align*}
\xi_{W} &= ((a_{X}^{\uparrow} \cdot a_{Y}^{\uparrow}) \cdot a_{Z}^{\uparrow}, (b_{X}^{\uparrow} \cdot b_{Y}^{\uparrow}) \cdot b_{Z}^{\uparrow}, (a_{X}^{\downarrow} \cdot a_{Y}^{\downarrow}) \cdot a_{Z}^{\downarrow}, (b_{X}^{\downarrow} \cdot b_{Y}^{\downarrow}) \cdot b_{Z}^{\downarrow}) \\
 &= (a_{X}^{\uparrow} \cdot (a_{Y}^{\uparrow} \cdot a_{Z}^{\uparrow}), b_{X}^{\uparrow} \cdot (b_{Y}^{\uparrow} \cdot b_{Z}^{\uparrow}), a_{X}^{\downarrow} \cdot (a_{Y}^{\downarrow} \cdot a_{Z}^{\downarrow}), b_{X}^{\downarrow} \cdot (b_{Y}^{\downarrow} \cdot b_{Z}^{\downarrow})) \\
 &= \xi_{U}
 \end{align*}
 The second equality results from the associativity of $\cdot$ on $\mathbb{R}$, hence, $W = U$.
 \item \textbf{Commutativity of Multiplication}: We have that 
 \begin{align*}
 \xi_{X} \boxtimes \xi_{Y} &= (a_{X}^{\uparrow} \cdot a_{Y}^{\uparrow}, b_{X}^{\uparrow} \cdot b_{Y}^{\uparrow}, a_{X}^{\downarrow} \cdot a_{Y}^{\downarrow}, b_{X}^{\downarrow} \cdot b_{Y}^{\downarrow}) \\
 &=  (a_{Y}^{\uparrow} \cdot a_{X}^{\uparrow}, b_{Y}^{\uparrow} \cdot b_{X}^{\uparrow}, a_{Y}^{\downarrow} \cdot a_{X}^{\downarrow}, b_{Y}^{\downarrow} \cdot b_{X}^{\downarrow}) \\
 &= \xi_{Y} \boxtimes\xi_{X}
 \end{align*}
\item \textbf{Multiplicative Identity:} The essential tuple $\xi_{1} = (1,1,1,1)$ is the multiplicative identity for $\mathbb{K}_{h,r}$, since multiplication is performed component-wise over $\mathbb{R}$.
\item \textbf{Distributivity of Multiplication with Respect to Addition:} 
We have that 
\begin{align*}
 \xi_X \boxtimes (\xi_Y \boxtimes \xi_Z) &= (a_{X}^{\uparrow} \cdot (a_{Y}^{\uparrow} + a_{Z}^{\uparrow}), b_{X}^{\uparrow} \cdot( b_{Y}^{\uparrow} + b_{Z}^{\uparrow}), a_{X}^{\downarrow} \cdot (a_{Y}^{\downarrow} + a_{Z}^{\downarrow}), b_{X}^{\downarrow} \cdot (b_{Y}^{\downarrow} + b_{Z}^{\downarrow})) \\
 &= (a_{X}^{\uparrow}a_{Y}^{\uparrow} + a_{X}^{\uparrow}a_{Z}^{\uparrow}, b_{X}^{\uparrow}b_{Y}^{\uparrow} + b_{X}^{\uparrow}b_{Z}^{\uparrow}, a_{X}^{\downarrow}a_{Y}^{\downarrow} + a_{X}^{\downarrow}a_{Z}^{\downarrow}, b_{X}^{\downarrow}b_{Y}^{\downarrow} + a{X}^{\downarrow}b_{Z}^{\downarrow}) \\
 &= (\xi_X \boxtimes \xi_Y) \boxplus (\xi_X \boxtimes \xi_Z)
\end{align*}
Right distributivity follows from commutativity of multiplication. 
\end{enumerate}
\end{proof}
\textit{Remark: Note that multiplicative inverses do not exist for all elements, and thus our stated division operation only applies to a subset of elements. This is why we do not have a field.}

We will illustrate this general construction for three specific sets of typed OFNs: trapezoidal, Gaussian, and exponential, as these three types are the most commonly used in applications such as~\cite{MB2013}.

\subsection{Trapezoidal OFNs}

For the base function $h(x) = x$, we call $\mathbb{K}_{l}$ the set of all true linear or true trapezoidal OFNs. Trapezoidal OFNs have up and down functions of the form $f(x) = ax + b$, $a,b \in \mathbb{R}$. 

OFNs in  $X \in \mathbb{K}_{l}$ have the form $X = (a_{X}^{\uparrow}x + b_{X}^{\uparrow}, a_{X}^{\downarrow} + b_{X}^{\downarrow})$. We give some examples of arithmetic on trapezoidal OFNs.

\begin{example}
Let $X = (-1,-x)$, $Y = (2,4x-2)$. Then $\xi_{X} = (0,-1,-1,0)$ and $\xi_{Y} = (0,2,4,-2)$.
We have that  $X \boxplus Y$ is the trapezoidal OFN $W$ with essential tuple $\xi_{W} = (0+0,-1+2,-1+4,0+(-2)) = (0,1,3,-2)$ and function representation $W = (1,3x-2)$. $X \boxdot Y$ is the trapezoidal OFN $Z$ with essential tuple $\xi_{Z} = (0\cdot 0, -1\cdot 2, -1\cdot 4, 0\cdot 2) = (1,-2,-4,0)$ and function representation $Z = (x-2,-4x)$. $X \squareslash Y$ is the trapezoidal OFN $U$ with essential tuple $\xi_{U} = (0/0,-1/2,-1/4,0/-2) = (0\cdot 0^{-1}, -1/2, -1/4, 0) = (1,-1/2,-1/4,0)$ with function representation $U = (x-1/2, -x/4)$.
\end{example}

\begin{example}
Let $X = (x-5,-x-3)$, $Y = (x+5, -x+3)$. Then $\xi_{X} = (1,-5,-1,-3)$, $\xi_{Y} = (1,5,-1,3)$. We have that $X \boxplus Y = (2x, -2x)$,  $X \boxminus Y = (-10, -6)$, $X \boxdot Y = (x - 25, x -9)$, and $X \squareslash Y = (x-1,x-1)$
\end{example}

 In Kosinski's original set of arithmetic operations, the product of two trapezoidal fuzzy numbers would not necessarily not result in a trapezoidal OFN. Our operations guarantee that the result of any operation on two trapezoidal OFNs remains a trapezoidal OFN. 
\subsection{A ring of Gaussian OFNs}
Let $h(x) = \sqrt{-2\log(x)}$. Then we call $\mathbb{K}_{G,r}$ the set of Gaussian OFNs if $X \in \mathbb{K}_{G,r}$ has the form $(a^{\uparrow}\sqrt{-2\log(x)} + b^{\uparrow}, a^{\downarrow}\sqrt{-2\log(x)} + b^{\downarrow})$. 

Again, the essential tuple for $X \in \mathbb{K}_{G,r}$ is $\xi_{X} = (a^{\uparrow}, b^{\uparrow}, a^{\downarrow}, b^{\downarrow})$, but this $\xi_{X}$ is defined with respect to $h(x) = \sqrt{-2\log(x)}$ rather than $h(x) = x$ as we did for the trapezoidal OFNs.

It should be noted that Gaussian OFNs do not have compact support, though we recognize instances where this may be desired, such as the use of OFNs for statistical applications.

\begin{figure}[H]
\begin{tikzpicture}

\begin{axis}[scale=0.9, xmin=-1, xmax=1.5, ymin=-.8, ymax=2.3, grid = both,
axis x line=middle,thick,xticklabels={},
axis y line=middle,tick style={draw=none},yticklabels={},
xlabel=$\alpha$, ylabel=$\mathbb{R}$,
every axis y label/.style={
    at={(ticklabel* cs:1.00)},
    anchor=west,
},
every axis x label/.style={
    at={(ticklabel* cs:1.00)},
    anchor=west,
}
]

\addplot[red, postaction={decorate, decoration={markings,
mark=at position 0.41 with {\arrow{>};},
mark=at position 0.71 with {\arrow{>};},
      }}] expression[domain=0:1, samples=100]{1/8*sqrt(-2*ln(x))+1};

\addplot[red, postaction={decorate, decoration={markings,
mark=at position 0.41 with {\arrow{<};},
mark=at position 0.71 with {\arrow{<};},
      }}] expression[domain=0:1, samples=100]{-1/8*sqrt(-2*ln(x))+1};
\node[below,text=red,font=\large] at (0.1,1.2){$Y$}; 

\addplot[domain=0:1,samples=100,blue,
postaction={decorate, decoration={markings,
mark=at position 0.41 with {\arrow{>};},
mark=at position 0.71 with {\arrow{>};},}}
        ]{1/4*sqrt(-2*ln(x))}node[right,pos=0.9]{};
        
\node[below,text=blue,font=\large] at (0.6,0.2){};
\addplot[domain=0:1,samples=100,blue,
postaction={decorate, decoration={markings,
mark=at position 0.31 with {\arrow{<};},
mark=at position 0.71 with {\arrow{<};},
    }}
        ]{-1/4*sqrt(-2*ln(x))}node[right,pos=0.9]{};

\addplot[domain=0:1,samples=100,violet,
postaction={decorate, decoration={markings,
mark=at position 0.41 with {\arrow{>};},
mark=at position 0.71 with {\arrow{>};},
      }}
        ]{3/8*sqrt(-2*ln(x))+1}node[right,pos=0.9]{};
\node[below,text=violet,font=\large] at (0.1,1.72){$Z$};
\addplot[domain=0:1,samples=100,violet,
postaction={decorate, decoration={markings,
mark=at position 0.31 with {\arrow{<};},
mark=at position 0.71 with {\arrow{<};},
    }}
        ]{-3/8*sqrt(-2*ln(x))+1}node[right,pos=0.9]{};
\node[below,text=blue,font=\large] at (.1,-.1){$X$};

\end{axis}
\end{tikzpicture}
\caption{\\ 
\color{blue}{$X = (\frac{1}{4}h(x), -\frac{1}{4}h(x))$},\\
\color{red}{$Y=(\frac{1}{8}h(x)+1, -\frac{1}{8}h(x)+1)$},\\
\color{violet}{$Z = X \boxplus Y = (\frac{3}{8}h(x)+1,
-\frac{3}{8}h(x)+1)$}}
\label{fig: add}
\end{figure}

\begin{example}
Let $X = (\frac{1}{4}h(x), -\frac{1}{4}h(x))$, 
$Y = (\frac{1}{8}h(x)+1, -\frac{1}{8}h(x)+1)$, where $h(x) = \sqrt{-2\log(x)}$.  Then $X \boxplus Y = (\frac{3}{8}h(x)+1,
-\frac{3}{8}h(x)+1)$ (Figure ~\ref{fig: add}).
\end{example}

Proper Gaussian OFNs may be inverted to yield the following membership functions. For increasing OFNs,
\[\mu_{A} = \begin{cases}
\exp\left(\frac{-(x-b^{\uparrow})^{2}}{2|a^{\uparrow}|}\right), & x \in (-\infty, b^{\uparrow}) \\
\exp\left(\frac{-(x-b^{\downarrow})^{2}}{2|a^{\downarrow}|}\right), & x \in (b^{\downarrow}, \infty) \\
1, & x \in [b^{\uparrow}, b^{\downarrow}]
\end{cases}\]

For decreasing OFNs, 
\[\mu_{A} = \begin{cases}
\exp\left(\frac{-(x-b^{\downarrow})^{2}}{2|a^{\uparrow}|}\right), & x \in (-\infty, b^{\downarrow}) \\
\exp\left(\frac{-(x-b^{\uparrow})^{2}}{2|a^{\uparrow}|}\right), & x \in (b^{\uparrow}, \infty) \\
1, & x \in [b^{\downarrow}, b^{\uparrow}]
\end{cases}\]

\subsection{The ring of exponential OFNs}

Let $h(x) = \log(x)$. Then the ring of exponential OFNs is denoted $\mathbb{K}_{e,r}$, and $X = (\mu_{A}^{\uparrow}, \mu_{A}^{\downarrow}) \in \mathbb{K}_{e,r}$ has the form $X = (a^{\uparrow}\log(x) + b^{\uparrow}, a^{\downarrow}\log(x) + b^{\downarrow})$ and essential tuple $\xi_{X} = (a^{\uparrow}, b^{\uparrow}, a^{\downarrow}, b^{\downarrow})$. 

\section{Improper Typed OFNs and Control of Pathologies}

In constructing rings of typed OFNs for given base functions $h$, we allow improper OFNs. In fact, in many cases, the multiplicative identity of a rings of typed OFNs is improper. Despite the presence of improper OFNs in a particular ring of typed OFNs, we remain at a practical advantage with our methods, since even the improper OFNs that exist in the ring are of the correct type. Moreover, typed OFNs allow for better identification and control of pathologies in the OFN that cause impropriety. 

Under Kosinski's operations, the product of two OFNs of a given type would not remain that same type. For example, the product of two trapezoidal OFNs under Kosinski's arithmetic operations would result in a quadratic OFN. Moreover, with such a broad definition of the space of OFNs and their operations, the resultant OFN could have up and down functions that are non-monotonic. 

Kacprzak and Kosinski~\cite{Kacprzak2011} do propose a method to "fix" improper OFNs so that these improper OFNs have a membership function and thereby become proper OFNs. For $A \in \mathbb{K}$ such that $A$ is improper, they define a new \textit{corresponding} membership function $\tilde{\mu}_{A}$ as 
\[\tilde{\mu}_{A}(x) = \begin{cases} \text{argmax}\{\mu_{A}^{\uparrow}(s) = x, \mu_{A}^{\downarrow}(y) = x\}, & x \in \text{Range}(\mu_{A}^{\uparrow}) \cup \text{Range}(\mu_{A}^{\downarrow}) \\
1, & x \in [\min\{\mu_{A}^{\uparrow}(1), \mu_{A}^{\downarrow}(1)\}, \max\{\mu_{A}^{\uparrow}(1), \mu_{A}^{\downarrow}(1)] \\
0, & \text{otherwise}
\end{cases}\]

While this is an elegant and compact fix for all pathologies, it does not preserve type, which we desire for application purposes. 

We revisit the properties an improper OFN may have, as described in Section~\ref{subsec: kosinski ofn}. Under our operations on a ring of typed OFNs where $h$ is always monotonic on the real line, such as trapezoidal, Gaussian, or exponential OFNs, the resultant OFN of any operation is always in the ring, and thus will never have a non-monotonic up or down function. Thus, (i) as a pathology is not possible under our operations, and we must only contend with fixing (ii) and (iii) as they arise. 

An OFN $A$ with essential tuple $\xi_{A} = (a^{\uparrow}, b^{\uparrow}, a^{\downarrow}, b^{\downarrow})$ as a pathology of type (ii) only if $\text{sgn}(a^{\uparrow}) = \text{sgn}(a^{\downarrow})$ and there is no $\alpha \in [0,1]$ such that $\mu_{A}^{\uparrow}(\alpha) = \mu_{A}^{\downarrow}(\alpha)$. In this case, we may correct the improper $A$ of type (ii) to a proper OFN $\hat{A}$ in the following way:

\begin{table}[H]
\begin{tabular}{|c|c|c|}
\hline 
Orientation & Sign & $A'$ \\
\hline
Increasing & $\text{sgn}(a^{\uparrow}) = \text{sgn}(a^{\downarrow}) = 1$ & $(\mu_{A}^{\uparrow}, b^{\downarrow})$ \\
Increasing & $\text{sgn}(a^{\uparrow}) = \text{sgn}(a^{\downarrow}) = -1$ & $(b^{\uparrow}, \mu_{A}^{\downarrow})$ \\
Decreasing & $\text{sgn}(a^{\uparrow}) = \text{sgn}(a^{\downarrow}) = 1$ & $(b^{\uparrow}, \mu_{A}^{\downarrow})$ \\
Decreasing & $\text{sgn}(a^{\uparrow}) = \text{sgn}(a^{\downarrow}) = -1$ & $(\mu_{A}^{\uparrow}, b^{\downarrow})$ \\
\hline
\end{tabular}
\end{table}

If $A$ is an improper OFN of type (iii) only, then $\text{sgn}(a^{\uparrow}) \neq \text{sgn}(a^{\downarrow})$, but there exists $\alpha \in [0,1)$ such that $\mu_{A}^{\uparrow}(\alpha) = \mu_{A}^{\downarrow}(\alpha)$. In this case, we correct improper OFNs $A$ of type (iii) to a proper OFN $\hat{A}$ given by 
\[\hat{A} = (a^{\uparrow}h + b^{\downarrow}, a^{\downarrow}h + b^{\uparrow}),\]
which "untwists" the pathology. 

If $A$ is a combination of types (ii) and (iii) (which will only happen on OFNs with compact support), then the correction $\hat{A}$ is given by $\hat{A} = (\mu_{A}^{\uparrow}(1), \mu_{A}^{\downarrow}(1)).$

\textit{Remark: In the applications where the result of an arithmetic operation of OFNs is required to be a classical fuzzy number, any needed corrections to the OFNs should be postponed until the end of all arithmetic operations in order to avoid propagation of error. One should perform arithmetic on the OFNs, be they proper or improper, and only correct the end result if the application requires it.}

\section{Conclusion and Future Work}

In this paper we have detailed the difficulties of current proposed operations on ordered fuzzy numbers, and noted that the algebraic structure induced by these various operations isn't as strong as applications would necessitate. We proposed the definition of typed OFNs, and detailed a method to construct a ring of typed OFNs for a given type. 

Improper OFNs can be problematic when an application requires the OFN to also have a membership function in the sense of a classical fuzzy number. A ring of typed OFNs does allow improper OFNs, but we gave a correction procedure for an improper OFN of a given type that forces the OFN to retain its type, unlike Kosinski's original correction in~\cite{Kacprzak2011}. 

While restricting the underlying set of a ring of typed OFNs to only one type may not yield the most elegant generalizations, it is extremely useful in applications. For instance, in~\cite{MB2013}, the authors utilized the concept of an OFN to create ordered fuzzy candlesticks from financial time series data, then proposed a generalization of the autoregressive model on the created OFNs. In practice, the user would choose a particular type of OFN to represent each window of data, and the autoregressive modelling would be performed on a set of typed OFNs (typically trapezoidal). Prokopowicz~\cite{Prokopowicz2016} applied OFNs to a fluid flow problem in a tank, focusing only on triangular OFNs (a subset of trapezoidal). Our operations allow both of these applications to perform arithmetic on these fixed types while guaranteeing the resultant remains of the desired type. It should be noted that ~\cite{Prokopowicz2016} provides one possible interpretation of improper OFNs, so correction would not be needed in this application. 
The theory of OFNs is still very sparse, and there remains much work
to be done.  One particular area of interest is ordered fuzzy
autoregression (OFAR), which requires a ring of OFNs to be 
endowed with a suitable metric; in fact, the comparison of
metrics on a space of OFNs is a question of interest in itself.

Another promising application lies in graph theory. Blue et al~\cite{Blue2002} describe generalizations of various path-finding and optimization algorithms commonly used on weighted graphs and networks that assign fuzzy weights to arcs. These problems can be expressed in terms of path algebras~\cite{Gondran1984} $(S, +, *)$ whose structure is generally only required to be a semi-ring. The set $S$ and the operations $+$ and $*$  vary based on the particular problem being solved. Common operations over crisp sets for $+$ and $*$ include minimum, maximum, union, conventional addition, and conventional multiplication. The minimum and maximum of two fuzzy numbers is taken pointwise on the endpoints of the level sets~\cite{KaufmanGupta1985}, and we have defined easily implemented addition and multiplication. Thus, we see application of our operations in implementing generalized graph algorithms with fuzzy arc weights. 

We plan to address these questions in future works.

\section{Acknowledgements}
Both authors would like to thank Jason Hathcock for many useful conversations and countless \LaTeX\ typesetting tips.  This research was not funded by any person nor organization; the authors declare no conflicts of interest.

\printbibliography
\end{document}